\documentclass[11pt]{amsart}
\usepackage{amssymb,amsfonts,amsmath}
\usepackage[all]{xy}

\newtheorem{thm}{Theorem}
\newtheorem{prop}[thm]{Proposition}
\newtheorem{conj}[thm]{Conjecture}
\newtheorem{qu}[thm]{Question}

\begin{document}

\title[Higher Chow groups]{A note on linear higher Chow groups}
\author{Muxi Li}
\address{Department of Mathematics, Washington University, St. Louis, MO  63130}
\email{muxili@math.wustl.edu}

\begin{abstract}
We give a counterexample to the proof in \cite{dJ02} of the existence of linear representatives of higher Chow groups of number fields.
\end{abstract}

\maketitle

\section{Introduction}
Higher Chow groups were introduced by Bloch \cite{Bl86} three decades ago to geometrize Quillen's higher algebraic $K$-theory. Let $X$ be a quasi-projective variety over an infinite field $k$. Writing $\Delta^m_k :=\mathbb{P}^m_k \setminus H\cong \mathbb{A}^m_k$, where $H:=\{x_0 + \cdots +x_m=0\}\cong \mathbb{P}^{m-1}$, denote by $Z^p(X,m)$ the free abelian group on closed irreducible subvarieties of $X\times \Delta^m_k$ of codimension $p$, properly intersecting each face $X\times \partial_I \Delta^m_k$ (given by $x_i =0$ $\forall i\in I$). One then defines $\mathrm{CH}^p(X,m)$ to be the $m^{\text{th}}$ homology of the complex $$\left( Z^p(X,\bullet),\,\partial:=\Sigma_{i=0}^{\bullet}(-1)^i \rho_i^*\right),$$where $\rho_i:\Delta^{\bullet-1}\hookrightarrow \Delta^{\bullet}$ is the inclusion of $\partial_i \Delta^{\bullet}$. If $X$ is smooth and $k$ is a subfield of $\mathbb{C}$, one has Bloch's Abel-Jacobi maps $$\mathrm{AJ}:\mathrm{CH}^p(X,m)\to H^{2p-m}_{\mathcal{H}}\left( X^{\text{an}}_{\mathbb{C}},\mathbb{Z}(p)\right)$$into absolute Hodge cohomology, which may be described ($\otimes \mathbb{Q}$) in terms of explicit maps of complexes $\widetilde{\mathrm{AJ}}$ \cite{BKLL15}. The homology of the subcomplex $\mathit{LZ}^p(X,\bullet)$ given by equations \emph{linear} in the $\{x_i\}$ defines the \emph{linear higher Chow groups} $\mathrm{LCH}^p(X,m)$, which map naturally to $\mathrm{CH}^p(X,m)$.

This note concerns the case $\mathrm{CH}^p(k,m)$ of a point over a number field, where $X=\text{Spec}(k)$. Working $\otimes\mathbb{Q}$, this is zero unless $(p,m)=(n,2n-1)$, in which case $\mathrm{CH}^n(k,2n-1)_{\mathbb{Q}}\cong K_{2n-1}(k)_{\mathbb{Q}} \cong K_{2n-1}(\mathcal{O}_k)_{\mathbb{Q}}$.  The linear group $\mathrm{LCH}^n(k,2n-1)_{\mathbb{Q}}$ is (for each $n\geq 1$) the image of a canonical homomorphism
$$\psi_n :\,H_{2n-1} \left( \mathrm{GL}_n(k),\mathbb{Q}\right)\to \mathrm{CH}^n(k,2n-1)_{\mathbb{Q}},$$
induced by the morphism of complexes
$$\tilde{\psi}_n:\, C_{\bullet}^{\text{grp}}(n)\to Z^n (k,\bullet)_{\mathbb{Q}}$$
given (for $\bullet = m$) by
$$(g_0,\ldots ,g_m) \longmapsto \left\{ \sum_{i=0}^m x_i \cdot g_i \underline{v} =0\right\} \subset \Delta^m$$
for some choice of $\underline{v}\in k^n \setminus \{\underline{0}\}$. (Here we consider $C_i^{\text{grp}}$ resp. $Z^n(k,i)$ to be in degree $-i$.) Now given an embedding $\sigma :k\hookrightarrow \mathbb{C}$, the Bloch-Beilinson regulator map (i.e., $\mathrm{AJ}$ composed with projection $\mathbb{C}/\mathbb{Q}(n)\twoheadrightarrow \mathbb{R}$) sends $\mathrm{CH}^n (\sigma(k),2n-1)_{\mathbb{Q}} \underset{r_{\text{Be}}}{\to} \mathbb{R}$, so that composing with \emph{all} $r=[k:\mathbb{Q}]=r_1 + 2r_2$ embeddings maps $\mathrm{CH}^n(k,2n-1)\to \mathbb{R}^r$. This factors through the invariants $\mathbb{R}^{d_n}$ [$d_n := r_2$ ($n$ even) resp. $r_1 + r_2$ ($n$ odd)] under de Rham conjugation, and is known to be equivalent to $\tfrac{1}{2}$ the Borel regulator $r_{\text{Bo}}: K_{2n-1}(\mathcal{O}_k)_{\mathbb{Q}}\to \mathbb{R}^{d_n}$ \cite{Bu02}.

Given the close relation between homology of $\mathrm{GL}_n$ and the original context of Borel's theorem, it is natural to consider the composite morphism of complexes $\widetilde{\mathrm{AJ}}\circ\tilde{\psi}_n$. Replacing $k$ by $\mathbb{C}$, these should yield explicit cocycles in $H_{\text{meas}}^{2n-1}\left( \mathrm{GL}_n(\mathbb{C}),\mathbb{C}/\mathbb{Z}(n)\right)$ lifting the Borel classes in $H^{2n-1}_{\text{cont}}\left( \mathrm{GL}_n(\mathbb{C}),\mathbb{R}\right)$ \cite{BKLL15}. This would also deepen our understanding of the equivalence of the Beilinson and Borel regulators. The first test of this proposal is to check its simplest implication:

\begin{conj}\label{thm1}
For a number field $k$, the linear higher Chow cycles surject {\rm (}rationally{\rm )} onto the simplicial higher Chow groups. Equivalently, $\psi_n$ is surjective for every $n\geq 1$.
\end{conj}

\section{A strategy for surjectivity?}
In fact, Conjecture \ref{thm1} is claimed as Proposition 16 in R. de Jeu's paper \cite{dJ02}. His approach is to fit (for each $n\geq 1$) $\tilde{\psi}_n$ into a commuting triangle
\begin{equation}\label{eq1}
\xymatrix{C_{\bullet}^{\mathrm{grp}}(n) \ar [r]^{\tilde{\psi}_n} \ar [rd]_{\tilde{r}_{\mathrm{Bor}}}  & Z^n(k,\bullet) \ar[d]^{\tilde{r}_{\mathrm{Be}}} \\ & \mathbb{R} [2n-1]. }
\end{equation}
Taking homology yields the diagram
\begin{equation}\label{eq2}
\xymatrix{H_{2n-1}(GL_n (k),\mathbb{Q}) \ar [r]^{\psi_n} \ar [rd]_{r_{\mathrm{Bor}}} & CH^n(k,2n-1)_{\mathbb{Q}} \ar [d]^{r_{\mathrm{Be}}} \\ & \mathbb{R} , }
\end{equation}
in which $r_{\mathrm{Bor}}$ [resp. $r_{\mathrm{Be}}$] is the Borel [resp. Beilinson] regulator, composed with a choice of embedding $k\hookrightarrow \mathbb{C}$. By composing with {\it all} embeddings (and using Borel's theorem), we get a diagram of the form
\begin{equation}\label{eq3}
\xymatrix{H_{2n-1}(GL_n(k),\mathbb{R}) \ar [r]^{\psi_n} \ar [rd]_{\cong} & CH^n(k,2n-1)_{\mathbb{R}} \ar [d]^{\cong} \\ & \mathbb{R}^{d_n} , }
\end{equation}
proving Conjecture \ref{thm1}.

The problem here is with de Jeu's choice of Goncharov's simplicial regulator $r_{\mathrm{Gon}}$ \cite{Go95} for $\tilde{r}_{\mathrm{Be}}$. While this appears to make \eqref{eq1} commute, by the calculation on pp. 228-230 of \cite{dJ02}, it is now known \cite{BKLL15} that $r_{\mathrm{Gon}}$ is not a map of complexes.  Specifically, in
\begin{equation}\label{eq4}
\xymatrix{\cdots \ar [r] & Z^2 (k,2n) \ar [d] \ar [r]^{\partial} & Z^n(k,2n-1) \ar [d]^{r_{\mathrm{Gon}}} \ar [r] & Z^n(k,2n-2) \ar [d] \ar [r] & \cdots \\ \cdots \ar [r] & 0 \ar [r] & \mathbb{R} \ar [r] & 0 \ar [r] & \cdots }
\end{equation}
we do not have $r_{\mathrm{Gon}} (\partial C_{2n}) = 0$.  So we must replace $r_{\mathrm{Gon}}$ by the ``corrected'' version in \cite{BKLL15}, which we will denote by $\mathrm{reg}_G$.  It is given on $Y\in Z^n(k,2n-1)$ by
\begin{equation}\label{eq5}
\mathrm{reg}_G (Y) := \int_{Y(\mathbb{C})} r_{2n-1} \left( \tfrac{x_1 + \cdots + x_{2n-1}}{-x_0}, \tfrac{x_2 + \cdots + x_{2n-1}}{-x_1},\ldots ,\tfrac{x_{2n-1}}{-x_{2n-2}} \right) ,
\end{equation}
which is known to induce $r_{\mathrm{Be}}$.

On the group homology side, de Jeu \cite{dJ02} also uses a formula of Goncharov for $\tilde{r}_{\mathrm{Bor}}$; we denote this by $\mathrm{reg}_B$.  Given $(g_0,\ldots ,g_{2n-1}) \in C^{\mathrm{grp}}_{2n-1}(n)$, let $\{ f_i \}_{i=1}^{2n-1}$ denote nonzero rational functions on $\mathbb{P}^{n-1}_{\mathbb{C}}$ with divisors
\begin{multline*}
D_i = \{ [\underline{X}]\in \mathbb{P}^{n-1} \mid (X_0,\ldots ,X_{n-1})\cdot g_i\underline{v} = 0\} \\
- \{ [\underline{X}]\in \mathbb{P}^{n-1} \mid (X_0,\ldots ,X_{n-1})\cdot g_0 \underline{v} = 0\} .
\end{multline*}
Then according to \cite[Thm. 5.12]{Go93},
\begin{equation}\label{eq6}
\mathrm{reg}_B (g_0,\ldots, g_{2n-1}) :=\int_{\mathbb{P}^{n-1}_{\mathbb{C}}} r_{2n-1} (f_1,\ldots ,f_{2n-1})
\end{equation}
induces $r_{\mathrm{Bor}}$. At least in the $n=2$ case we treat below, this formula is correct.   (See the calculation in $\S$\ref{s3} below.)  Moreover, it is well-defined for any $n$, in the sense that the RHS of \eqref{eq6} is invariant when we rescale any $f_i$ by a constant.

We tried to emulate the approach in \cite{dJ02} to see if the new diagram \eqref{eq1} (with $\tilde{r}_{\mathrm{Bor}} = \mathrm{reg}_B$ unchanged and $\tilde{r}_{\mathrm{Be}}$ corrected to $\mathrm{reg}_G$) commutes, with no success. At this point, we decided to attempt the first nontrivial case by hand, and arrived at a negative result:

\begin{prop}\label{pr1}
For $n=2$, the amended triangle \eqref{eq1} does not commute.
\end{prop}

\section{Proof of Proposition \ref{pr1}}\label{s3}
In \cite{Go04}, Goncharov mentions the formula
\begin{equation}\label{eq7}
\int_{\mathbb{P}^1} r_3 (f_1,f_2,f_3) = \sum_{(x_1,x_2,x_3) \in \mathbb{C}^3} \nu_{x_1}(f_1)\nu_{x_2}(f_2)\nu_{x_3}(f_3) D_2(\mathsf{CR}(x_1,x_2,x_3,\infty))
\end{equation}
where $\nu_x(f)$ is the order of $f$ at $x$.  One easily verifies that this is correct; it will be required below.

Now take $\underline{v}=\left( \begin{array}{cc} 1 \\ 0 \end{array} \right) \in \mathbb{C}^2$, and $(g_0,g_1,g_2,g_3) \in C_3 (2)$. We can do a change of coordinate to let $g_0=\left( \begin{array}{cc} 1 & * \\ 0 & * \\ \end{array} \right)$, $g_1=\left( \begin{array}{cc} 0 & * \\ 1 & * \\ \end{array} \right)$, $g_2=\left( \begin{array}{cc} a & * \\ c & * \\ \end{array} \right)$, $g_3=\left( \begin{array}{cc} b & * \\ d & * \\ \end{array} \right)$. For convenience, we set $\Delta:=ad-bc$.

Write $z:=\tfrac{X_1}{X_0}$ and $f_1(z)=z$, $f_2(z)=cz+a$, and $f_3(z)=dz+b$.  According to \eqref{eq6} and \eqref{eq7}, we have
\begin{flalign*}
\mathrm{reg}_B (g_0,g_1,g_2,g_3) &= \int_{\mathbb{P}^1} r_3 (z,cz+a,dz+b) \\
&= D_2 \left( \tfrac{bc}{ad} \right) .
\end{flalign*}
This is consistent with evaluating the cocycle $\varepsilon_2 \in H^3_{\mathrm{cont}}(GL_2(\mathbb{C}),\mathbb{R})$ (cf. Intro. to \cite{BKLL15}) on the ``group homology chain'' $(g_0,g_1,g_2,g_3)$.

For the other side, applying $\tilde{\psi}$ to this chain produces the linear higher Chow chain $Y \subset \Delta^3$ cut out by
\[  x_0 + ax_2 +bx_3 = 0 \;\;\text{and} \;\; x_1 + cx_2 + dx_3 = 0. \]
Parametrizing $Y\cong \mathbb{P}^1$ by $t\mapsto (\Delta ,\Delta t ,bt-d ,c-at)$, \eqref{eq5}, \eqref{eq7} and the rescaling property yield $\mathrm{reg}_G (\tilde{\psi}(g_0,g_1,g_2,g_3) ) =$
\begin{flalign*}
\mathrm{reg}_G (Y) &= \int_{\mathbb{P}^1} r_3 \left( \tfrac{(d-c)+(a-b-\Delta)t}{\Delta}, \tfrac{(d-c)+(a-b)t}{\Delta t},\tfrac{at-c}{bt-d}\right) \\
&= \int_{\mathbb{P}^1} r_3\left( (c-d)+(\Delta + b-a)t, \tfrac{(c-d)+(b-a)t}{t}, \tfrac{at-c}{bt-d}\right) \\
&= D_2\left(\tfrac{(d-1)\Delta}{b(c-d)}\right) - D_2\left(\tfrac{(c-1)\Delta}{a(c-d)}\right) -D_2\left(\tfrac{(b-a)(d-1)}{b(d-c)}\right) + D_2\left( \tfrac{(b-a)(c-1)}{a(d-c)}\right) .
\end{flalign*}

To check that these two results disagree, put $a=1$, $b=-1$, $c=1-i$, $d=1+i$, so that $\Delta = 2$ and $\tfrac{ad}{bc}=-i$.  Of course, $D_2(-i)\neq 0$.  On the other hand,
\[ \tfrac{(d-1)\Delta}{b(c-d)} , \; \tfrac{(c-1)\Delta}{a(c-d)}, \; \tfrac{(b-a)(d-1)}{b(d-c)},\; \tfrac{(b-a)(c-1)}{a(d-c)} \]
are all $1$, $D_2$ of which is $0$.

\section{Concluding remarks}
Naturally, it is still possible that \eqref{eq2} commutes, since there we restrict to \emph{closed} chains.  In fact, even if we don't accept the proof in \cite{dJ02}, there is the earlier result of Gerdes \cite{Ge91} which gives surjectivity of $\psi_n$ for $n=2$.  Moreover, there is the agreement between the Beilinson and Borel regulators in \cite{Bu02}, though this does not involve $\psi_n$ in any way.  To sum up, we conclude with the
\begin{qu}
Are there any techniques to prove that \eqref{eq2} commutes even though the amended diagram \eqref{eq1} does not, for $n=2$ and more generally?  Or is it more likely that $\psi_n$ has to be somehow modified?
\end{qu}

\subsection*{Acknowledgments}
This work was supported by the National Science Foundation [DMS-1361147; PI: Matt Kerr].  The author would like to thank his advisor Matt Kerr and Jos\'e Burgos Gil for discussions regarding this note.

\end{document}